\documentclass[11pt]{article}
\usepackage{amssymb}
\usepackage{amsmath}
\usepackage{amsthm}

 \usepackage{color}
 \usepackage{ulem}

\numberwithin{equation}{section}

\date{}

\title{\bf  Notes on nonlocal dispersal equations in  a periodic habitat}

\author{Jian Fang
\\Department of Mathematics\\
 Harbin Institute of Technology\\
  Harbin 150001, China\\
  E-mail:\, jfang@hit.edu.cn
\and
Xiao-Qiang Zhao\\
Department of Mathematics and Statistics\\
Memorial University of Newfoundland\\
St. John's, NL A1C  5S7, Canada\\
E-mail:\,  zhao@mun.ca}

\begin{document}
\maketitle


Assume that $\int_{\mathbb{R}}J(y)dy=1$.  We consider the linear
equation
\begin{equation}\label{E1}
u_t=J*u-u, \quad  t>0, \, x\in \mathbb{R}.
\end{equation}
Let $b>0$ be given, and $\mathcal{C}_b:=C(\mathbb{R},[0,b])$.
For any $\phi, \psi\in \mathcal{C}_b$, we define 
$$
d(\phi,\psi)=\sum_{k=1}^{\infty}\frac 1{2^k}\max_{-k\leq x\leq k}|\phi(x)-\psi(x)|.
$$
Then  $(\mathcal{C}_b, d)$ is a complete metric space, and the metric $d$
induces the compact open topology on $\mathcal{C}_b$. Let $T(t):\mathcal{C}_b\to
\mathcal{C}_b$ be the solution maps associated with \eqref{E1}.

\

\noindent {\bf Lemma 1.} For each $t>0$, the map $T(t):
\mathcal{C}_b\to \mathcal{C}_b$ is an $\alpha$-contraction,
where $\alpha$ is the Kuratowski measure of noncompactness in the metric
space $(\mathcal{C}_b, d)$.

\noindent {\bf Proof.}  Letting $n=1$ in (3.2) and (3.5) of
\cite{WZ} with $P(x)=J(x)$, we see that
\begin{equation}
T(t)\phi(x)=e^{-t}\sum_{k=0}^{\infty}\frac{t^k}{k!}a_k(\phi)(x),
\quad \forall t\geq 0, \, x\in \mathbb{R},
\end{equation}
where
$$
a_0(\phi)(x)=\phi(x), \quad
a_k(\phi)(x)=\int_{\mathbb{R}}J(x-y)a_{k-1}(\phi)(y)dy,\, \forall
k\geq 1.
$$
By induction, we have  $a_k(\phi)(x)\leq b, \, \forall x\in \mathbb{R}, \, \phi\in \mathcal{C}_b$.
Define
$$g(x)=\int_{\mathbb{R}}|J(x+z)-J(z)|dz.$$
Then we obtain
$$
|a_k(\phi)(x_1)-a_k(\phi)(x_2)|\leq bg(x_1-x_2), \forall x_1, x_2\in \mathbb{R}, \, \phi\in \mathcal{C}_b,
\, k\geq 1.
$$
Let
$$T_1(t)\phi(x)=e^{-t}\phi(x),\quad
T_2(t)\phi(x)=e^{-t}\sum_{k=1}^{\infty}\frac{t^k}{k!}a_k(\phi)(x).
$$
Clearly, $T(t)=T_1(t)+T_2(t)$. Since
$e^t=\sum_{k=0}^{\infty}\frac{t^k}{k!}$, we have
$$
|T_2(t)\phi(x_1)-T_2(t)\phi(x_2)|\leq e^{-t}
bg(x_1-x_2)\sum_{k=1}^{\infty}\frac{t^k}{k!}\\
\leq bg(x_1-x_2),
$$
for all $x_1, x_2\in \mathbb{R}, \, \phi\in \mathcal{C}_b$. By
Ascoli's theorem, it then easily follows that for each $t>0$, the
map $T_2(t)$ is compact with respect to the compact open topology.
Note that $T_1(t)$ is an $\alpha$-contraction with the contraction
coefficient being $e^{-t}$. Consequently, $T(t)$ is also an
$\alpha$-contraction with the contraction coefficient being
$e^{-t}$.

\

From the proof of Lemma 1, we have the following observation.

\

\noindent {\bf Lemma 2.}  For any given nonempty and bounded interval $I:=[a,b]\subset\mathbb{R}$ and $t\geq 0$, there holds
$\alpha ((T(t)\mathcal{U})_I)\leq e^{-t}\alpha  ((\mathcal{U})_I)$
for each set $\mathcal{U}\subset \mathcal{C}_b$, 
where $\alpha$ is the Kuratowski measure of noncompactness in the Banach 
space $C([a,b],\mathbb{R})$.

\

By the basic properties of the  noncompactness measure, it is easy to prove 
the following result (see, e.g.,  the proof of \cite[Lemma 5]{Banas1980}).
 
 \
 
\noindent
{\bf Lemma 3.}  Let $a<b$ be two real numbers and $X$ be a Banach space. 
For a bounded set $B\subset C([a,b],X)$,  define 
$B(s):=\{g(s):  g\in B\}, \, s\in [a, b]$,  and $\int_a^b B(s){\rm d}s:=\{\int_a^b g(s){\rm d}s: g\in B\}$.  If $B$ is equicontinuous  on $[a,b]$, 
then $\alpha\left(\int_a^b B(s){\rm d}s\right)\leq \int_a^b \alpha(B(s)){\rm d}s$.

\

Next we consider the nonlinear equation
\begin{equation}\label{nonL}
u_t=J*u-u +f(x,u) \quad  t>0, \, x\in \mathbb{R}.
\end{equation}
Assume that $f(x,0)\equiv 0$ and $f(x+L,u)=f(x,u)$ for some $L>0$, system \eqref{nonL} has a positive $L$-periodic steady state $\beta(x)$, and there exists a real number $k_f>0$ such  that 
\begin{equation}
\label{Lip}
|f(x,u_1)-f(x,u_2)|\leq k_f|u_1-u_2|, \,  
\forall x\in [0,L],  u_1, u_2\in [0,\beta (x)].
\end{equation}

Let $Q(t)$ be the solution semiflow of
system \eqref{nonL} on $\mathcal{C}_{\beta}$.  Then we have the 
following result.

\

\noindent {\bf Theorem 1.}   For any given  nonempty interval $I=[a,b]\subset \mathbb{R}$ and  $t\geq 0$, there holds
$\alpha ((Q(t)B)_I)\leq e^{(k_f-1)t}\alpha  (B_I)$
for each set $B\subset \mathcal{C}_{\beta}$, 
where $\alpha$ is the Kuratowski measure of noncompactness in the Banach 
space $C([a,b],\mathbb{R})$.

\noindent{\bf Proof.}  Define $\hat f (\phi)(x)=f(x,\phi(x)), \, \forall x\in \mathbb{R}$.
 By the constant-variation formula, it follows that
\[
	Q(t)\phi=T(t)\phi +\int_0^tT(t-s)\hat f(Q(s)\phi)ds, \, \, \forall t>0.
\]
Define $Q_t\phi=Q(t)\phi$ and  let $B\subset \mathcal{C}_{\beta}$ be given.
In view of  Lemmas 2 and 3, we then have
\begin{eqnarray*}
	\alpha( (Q_tB)_I)&&\le \alpha ((T(t)B)_I)+\int_0^t\alpha((T(t-s)\hat  f(Q_sB))_I)ds\\
	&& \le  e^{- t}\alpha(B_I)+\int_0^t  e^{-(t-s)}\alpha((\hat f(Q_sB))_I)ds\\
	&& \le e^{- t}\alpha(B_I)+\int_0^t e^{-(t-s)}k_f\alpha((Q_sB)_I)ds,
\end{eqnarray*}
and hence,
\[
e^{t} \alpha((Q_tB)_I)= \alpha(B_I)+\int_0^t k_f e^{ s}\alpha((Q_sB)_I)ds.
\]
By Grownwall's inequality,  it follows that
\[
e^{t} \alpha((Q_tB)_I)\le  \alpha(B_I)e^{\int_0^t k_f ds}=\alpha(B_I)e^{k_f t}.
\]
This implies that  $\alpha(Q_t(B)_I)\le  e^{(k_f-1)t}\alpha(B_I)$.

\

\noindent
{\bf Remark 1.} From the proof of Theorem 1, it is easy  to see that the conclusion 
of Theorem 1 also holds if the condition \eqref{Lip} is replaced by the 
assumption that $\alpha((\hat f(B))_I)\leq k_f\alpha((B)_I)$ for any $B\subset \mathcal{C}_{\beta}$.

\

\noindent
{\bf Remark 2.}  Let  $Q_t$ be the solution maps of $u_t=D(J\ast u-u)+f(x,u)$.
By the proof of Theorem 1, we see that for each $t>0$,  $Q_t$ is an $\alpha$-contraction relative to $I$ provided that $D> k_f$. It follows that 
the large dispersal rate implies the $\alpha$-contraction property of solution 
maps for this class of nonlinear nonlocal dispersal equations.

\

\noindent 
{\bf Remark 3.} The spreading speed and traveling waves were 
obtained in \cite{ShenZhangJDE2010,ShenZhangCANA2012,CovPoincare2013,RSZDCDSA2015}
for nonlocal dispersal equation \eqref{nonL} under the condition that $\frac{f(x,u)}{u}$ is strictly
decreasing in $u>0$, which guarantees the linear determinacy of the spreading speed. 
In the case where $k_f<1$,  Theorem 1 implies that the solution
map $Q_t$ satisfies the  $\alpha$-contraction assumption (E3) in \cite{LZ} 
for each $t>0$. Thus,  one can use \cite[Theorems 5.2 and 5.3]{LZ}
to  establish the existence of  the spreading speed and their coincidence
with the minimum wave speed for monotone  spatially $L$-periodic traveling 
waves of nonlocal  equation \eqref{nonL} no matter whether
the spreading speed is linearly  determinate. 

\

\noindent 
{\bf Remark 4.}  It is easy to generalize  Theorem 1 and Remark 2  to
 nonlocal dispersal systems in  a periodic habitat.  By  the results and arguments in \cite{LZ,FZ2014},  it then follows that  the nonlocal dispersal 
 Lotka-Volterra competition system, as studied in  \cite{KRS2015, BaoLiShen2016}, admits a spatially periodic traveling 
 wave with the wave speed being the spreading speed under the condition that the dispersal rates are large.  However, it remains an 
 open problem whether this system admits such a critical wave without assuming the large dispersal rates.

\end{document}